\documentclass[12pt]{amsart}

\usepackage{amsmath,amsfonts,amsthm,amsopn,cite,mathrsfs}
\usepackage{epsfig,verbatim}
\usepackage{subfigure}

\newcommand{\tfa}{time-frequency analysis}

\newcommand{\fif}{if and only if}

\newtheorem{tm}{Theorem}

\newtheorem{cor}[tm]{Corollary}

\theoremstyle{definition}

\newcommand{\beqa}{\begin{eqnarray*}}
\newcommand{\eeqa}{\end{eqnarray*}}

\newcommand{\field}[1]{\mathbb{#1}}
\newcommand{\bR}{\field{R}}        
\newcommand{\bN}{\field{N}}        
\newcommand{\bZ}{\field{Z}}        
\newcommand{\bC}{\field{C}}        
\def\cO{\mathcal{ O}}

\def\<{\left<}
\def\>{\right>}

\def\inv{^{-1}}

\def\mv1{M_v^1}

\newcommand{\tp}{totally positive}
\newcommand{\lpc}{Laguerre-Polya class}
\newcommand{\pff}{Polya frequency function}

\begin{document}
\begin{abstract}
  We review Schoenberg's characterization of totally positive
  functions and its connection to the Laguerre-Polya class. This characterization
yields a new condition that is equivalent to the truth of the Riemann hypothesis.
\end{abstract}

\title[Totally positive functions and the zeta function]{
  Schoenberg's Theory of 
  Totally Positive Functions and the Riemann Zeta Function} 
\author{Karlheinz Gr\"ochenig}
\address{Faculty of Mathematics \\
University of Vienna \\
Oskar-Morgenstern-Platz 1 \\
A-1090 Vienna, Austria}
\email{karlheinz.groechenig@univie.ac.at}
\subjclass[2010]{}
\date{}
\keywords{}
\thanks{K.\ G.\ was
  supported in part by the  project P31887-N32  of the
Austrian Science Fund (FWF)}
\maketitle

In a series of papers in the 1950s Schoenberg investigated the
properties of \tp\ functions~\cite{CS66,sch47,Sch50,sch51,SW53}. He found several
characterizations and used total positivity to prove fundamental
properties of splines~\cite{SW53,Sch73}. 
The purpose of this note
is to survey some aspects of Schoenberg's work on \tp\ functions, 
to  advertize the connection between \tp\ functions and the Riemann
hypothesis,  and to provide some mathematical entertainment. 

One may speculate whether Schoenberg himself thought
about the Riemann zeta function.  He was the son-in-law of the
eminent number theorist Edmund Landau, he  collaborated  with Polya, he  knew deeply the work of Polya
and Schur about the \lpc\ of entire functions that remains influential
in the study of the Riemann hypothesis.  Yet, to
my knowledge, he  never mentioned any number theory in his work on
\tp\ functions and splines; by the same token,  Schoenberg's name is
not mentioned in analytic number theory. 

\vspace{3mm}

\emph{Totally positive functions.} A measurable function $\Lambda $ on $\bR $  is \tp , if for every
$n\in \bN$ and every two sets of
increasing numbers $x_1 < x_2 <\dots < x_n$ and $y_1< y_2 < \dots <
y_n$ the matrix $\Big( \Lambda (x_j - y_k) \Big) _{j,k=1, \dots ,n}$
has non-negative determinant:
\begin{equation}
  \label{eq:1}
  \det ( \Lambda (x_j - y_k) \Big) _{j,k=1, \dots ,n} \geq 0 \, .
\end{equation}
If in addition $\Lambda $ is integrable, then $\Lambda $
is called a Polya frequency function.

If  $\Lambda $ is \tp\ and not equal to $e^{ax+b}$, there exist an exponential
$e^{cx}$, such that $\Lambda _1(x) = e^{cx}\Lambda (x)$ is a \pff ,
i.e., $\Lambda _1$ is \tp\ and integrable~\cite[Lemma~4]{sch51}. It is
usually  no loss of generality to restrict to \pff s.

The class of \tp\ functions played and plays an important role in
approximation theory, in particular in spline theory~\cite{SW53}, and in
statistics~\cite{Efr65,karlin}. In a different and rather surprising
direction, \tp\ functions appear in the representation theory of
infinite dimensional motion groups~\cite{Pick91}. Recently, \tp\
functions appeared in sampling theory and 
in \tfa\ \cite{GS13,GRS18,GRS20}, where they were instrumental in the 
derivation of  optimal results.

\vspace{3mm}

\emph{The Laguerre-Polya class.} An entire function $\Psi $ of order
at most $2$ belongs to the
Laguerre-Polya class, if its Hadamard
factorization   is of the form 
\begin{equation}
  \label{eq:2}
  \Psi (s) = C s^m e^{-\gamma s^2 + \delta s} \, \prod _{j=1}^\infty
  (1+\delta _j s) e^{-\delta _js} \qquad s\in \bC \, ,
\end{equation}
where $\delta _j\inv \in \bR $ are the zeros of $\Psi $, $m$ is the
order of the zero at $0$, $\gamma \geq 0$, $\delta \in \bR $, and
\begin{equation}
  \label{eq:3}
  0<\gamma + \sum _{j=1}^\infty \delta _j^2 <\infty \, .
\end{equation}
Thus the \lpc\ consists of entire functions of order two with
convergence exponent two with only real zeros. 
While the study of the distribution of zeros of entire functions is a
perennial topic in complex analysis and of interest in its own right~\cite{Levin80},
the \lpc\ has gained special prominence in analytic number theory:  the Riemann
hypothesis says that a relative of the
Riemann zeta function belongs to the \lpc .

\vspace{3mm}

\emph{Schoenberg's characterization of \tp\ functions.}
The fundamental results about \tp\ functions were derived by
Schoenberg in  a series of papers~\cite{sch47,Sch50,sch51,SW53}. A
comprehensive treatment is contained in  Karlin's  monograph~\cite[Ch.~7]{karlin}. 

The notions of
\tp\ functions and \lpc\ are 
seemingly unrelated, yet there is a deep connection between them
through 
the following characterization of  Schoenberg~\cite{sch51}. 

\begin{tm} \label{tm:tp}
  (i) If $\Lambda $ is a Polya frequency function, then its
  (two-sided) Laplace transform converges in a vertical strip $\{z\in
  \bC : \alpha < \mathrm{Re}\, z < \beta \}, \alpha <0< \beta $, and
  \begin{equation}
    \label{eq:5}
    \int _{-\infty }^\infty \Lambda (x) e^{-sx} \, dx = \frac{1}{\Psi (s)} 
  \end{equation}
  is the reciprocal of a function $\Psi $ in the \lpc\ with $\Psi (0)>0$.

  (ii) Conversely, if $\Psi $ is in the \lpc\ with $\Psi (0) >0$, then
  its reciprocal 
  $1/\Psi $ is the Laplace transform of a Polya frequency function
  $\Lambda $.   
\end{tm}
This is a fascinating theorem,
because it relates two function classes that seem to bear absolutely no
resemblance to each other.  Schoenberg's theorem establishes a bijection between
the class of \pff s, the \lpc , and yields a  parametrization by the set $(0,\infty )
\times \bR \times \ell ^2(\bZ )$. 

By using the Fourier transform instead of the Laplace transform,
Schoenberg's theorem 
can  be recast as follows:
A function $\Lambda $ is \tp\ and integrable, \fif\ its Fourier
transform possesses the factorization
\begin{equation}
  \label{eq:6}
  \hat{\Lambda}(\tau  ) = C e^{-\gamma \tau  ^2 + 2\pi i \delta \tau  }\,
  \prod _{j=1}^\infty (1+2\pi i \delta _j\tau  )\inv e^{2\pi i \delta _j
    \tau  }
\end{equation}
where $C>0$, $\gamma \geq 0$, $\delta , \delta _j \in \bR$ and $\sum
_{j=1}^\infty \delta _j^2 < \infty $ (and the product in \eqref{eq:6}  may also be finite).  

If we drop the condition of integrability and exclude exponential
functions, then the representation \eqref{eq:6} still
holds for every \tp\ function, but their  Laplace transform of $\Lambda $ 
converges in some   vertical strip $\{z\in
  \bC : \alpha < \mathrm{Re}\, z < \beta \}$ that does not contain
  $0$.

  A similar result holds for one-sided \tp\ functions~\cite[Thm.~2]{sch51}.  
  \begin{tm} \label{tm:onesided}
(i) If $\Lambda $ is a Polya frequency function with support in
$[0,\infty )$, then its
   Laplace transform converges in a half-plane  $\{z\in
  \bC : -\alpha < \mathrm{Re}\, z  \}, \alpha  >0$, and
  \begin{equation}
    \label{eq:5}
    \int _0^\infty  \Lambda (x) e^{-sx} \, dx = \frac{1}{\Psi (s)} 
  \end{equation}
  is the reciprocal of an entire function $\Psi $ with Hadamard factorization
  \begin{equation}
    \label{eq:ne1}
  \Psi (s) = Ce^{\delta s} \prod _{j=1}^\infty (1+\delta _js) \, ,  
  \end{equation}
  with $\delta \in \bR , \delta _j \geq 0, \sum \delta _j <\infty $.
  
  (ii) Conversely, if $\Psi $ possesses the factorization \eqref{eq:ne1}, then
  its reciprocal 
  $1/\Psi $ is the Laplace transform of a Polya frequency function
  $\Lambda $ with support in $[0,\infty )$.       
  \end{tm}

  \vspace{3mm}
  
\emph{Elementary examples.}  If $\hat{\Lambda }(\tau  ) = (1+2\pi i \delta \tau  )\inv $,
then $\Lambda (x) = \delta \inv e^{-x/\delta }  \chi _{[0,\infty
  )}(x)$ is the one-sided exponential function. For $\hat{
  \Lambda}(\tau  ) =
e^{-\pi \gamma \tau  ^2}$ for $\gamma >0$, we obtain the
Gaussian $\Lambda (x) = \gamma ^{-1/2} e^{-\pi x^2/\gamma }$. In both
cases, it is easy to check directly that these functions are \tp .

\vspace{3mm}

The
proof of the implication  (ii) of Theorem~\ref{tm:tp}   is based on
the (non-trivial) fact that the convolution $\Lambda = 
\Lambda _1 \ast \Lambda _2$ of two \pff s $\Lambda _1 , \Lambda _2$ is
again a \pff .
The converse in Theorem~\ref{tm:tp}  lies much deeper,  and  Schoenberg
used heavily several results of Polya 
about functions in the \lpc ~\cite{Pol15,PS1914}. See the end of this
note for the essential step of the argument.  

Schoenberg's motivation  was  the characterization and deeper
understanding of \tp\ functions, and thus  the implication (i) and the
factorization \eqref{eq:6} can be considered  his  main insight about
\tp\ functions. However,  instead of reading Schoenberg's theorem as a
characterization of \tp\ functions, one may read it as a
characterization of the \lpc . \emph{A
  function $\Psi$ with $\Psi (0)>0$ and $\Psi \neq e^{as+b}$ is in the \lpc , \fif\ the Fourier transform of
  $1/\Psi $ is a \pff .}  

\vspace{3mm}
\emph{The Riemann hypothesis and \tp\ functions.}
Let $\zeta (s) = \sum _{n=1}^\infty n^{-s}$ for $s\in \bC ,
\mathrm{Re}\, s >1,$ be the Riemann zeta function and let
\begin{align}
    \xi (s) &= \tfrac{1}{2} s(s-1) \pi ^{-s/2} \Gamma
              \Big(\frac{s}{2}\Big) \, \zeta (s)   \label{eq:6b} \\
  \Xi (s) &= \xi \Big(\frac{1}{2} + is \Big) \label{eq:7}
\end{align}
be the Riemann xi-functions (where $\Gamma $ is the usual gamma
function). Then the functional equation for the Riemann zeta function
is expressed by the symmetry
\begin{equation}
  \label{eq:8}
    \xi (s) = \xi (1-s) \qquad \text{ and } \qquad \Xi (s) = \Xi (-s)
  \end{equation}
  for the xi-functions.
  
The Riemann hypothesis conjectures that  all 
non-trivial zeros of the zeta function lie  on the critical line $1/2
+ i t$. See the monographs~\cite{Iwa14,Ivi03,Tit86}, 
the two volumes about equivalents of the Riemann hypothesis~\cite{Br17-1,Br17-2} or the  survey articles~\cite{Bom10,Conrey03}.

Expressed in terms of  the
xi-functions, the Riemann hypothesis states that $\Xi $ has only real
zeros, in other words, \emph{$\Xi $ belongs to the \lpc .} Thus many
investigations of the zeta function involve complex analysis 
related to the \lpc .  Schoenberg's theorem immediately leads to  the following equivalent
condition for the Riemann hypothesis to hold.

\begin{tm} \label{tm:equi1}
  The Riemann hypothesis holds, \fif\ there exists a \pff\ $\Lambda $, such that 
  \begin{equation}
    \label{eq:9}
    \frac{1}{\Xi (s) } = \int _{-\infty } ^\infty \Lambda (x) e^{-sx} 
    \, dx \qquad \text{ for } s\in \bC , |\mathrm{Re} \, s| < t_0   \, , 
  \end{equation}
  where $1/2+i t_0$ is the first zero of the zeta function on the
  critical line. 
\end{tm}

Let us make this statement a bit more explicit by taking the Fourier
transform instead of the Laplace transform.

\begin{tm} \label{tm:equi2}
The Riemann hypothesis holds, \fif\ 
  \begin{equation}
    \label{eq:10}
    \Lambda (x) =  \frac{1}{2\pi}\int _{-\infty } ^\infty \frac{1}{\xi
      (\tfrac{1}{2}+ \tau  )
    } e^{- i x\tau } \, d\tau
  \end{equation}
is a Polya frequency  function.  
\end{tm}
  The growth of $\xi $ in the complex plane is~\cite{Tit86}
  $
  |\xi (s) | = \cO \Big( e^{A |s| \ln |s|}\Big) \, ,
  $
  and on the positive  real line
  $$\ln \xi (\sigma )  \asymp \tfrac{1}{2} \sigma  \log \sigma  \qquad
  \sigma >1 \, .
  $$
  Consequently $\frac{1}{\xi(\sigma )} \leq C e^{-|\sigma | \log
    |\sigma | /2}$ decays super-exponentially. Since  $\zeta $
  and thus $\xi $ do not have any real zeros in the interval
  $[0,1]$ and $\zeta >0$ on $(1,\infty)$, the function $\xi $ is therefore strictly positive on $\bR $ and
$1/\xi $ is   integrable. Thus its   Fourier transform is
well-defined. 
  
 Using $s=2\pi i \tau $, we can rewrite  \eqref{eq:9} as a 
 Fourier transform.  
 The inversion formula
  for the Fourier transform now yields
  \begin{align*}
        \Lambda (x) &=  \int _{-\infty } ^\infty \frac{1}{\Xi (2\pi i \tau  )}
                      e^{2\pi  i x\tau } \, d\tau  \\
    &= \int _{-\infty } ^\infty \frac{1}{\xi (1/2 - 2\pi  \tau  )}
                      e^{2\pi  i x\tau } \, d\tau \, ,
   \end{align*}
which is   \eqref{eq:10}. 

\vspace{3mm}

Using the symmetry of $\Xi $, there is an alternative formulation of
Theorem~\ref{tm:equi1} with the restricted \lpc\ defined in \eqref{eq:ne1}. Since
$\Xi $ is symmetric, it can be written as $\Xi (s) = \Xi _1(-s^2)$ for
an entire function $\Xi _1$ of order $1/2$. Furthermore, $\Xi $ has
only real zeros, \fif\ $\Xi _1$ has only negative zeros (with
convergence exponent at most $1$). The characterization of one-sided
\pff s yields the following equivalence.

\begin{tm} \label{tm:ones}
The Riemann hypothesis holds, \fif\ there exists a \pff\ $\Lambda $
with support in $[0,\infty )$, such that 
  \begin{equation}
    \label{eq:9}
    \frac{1}{\Xi _1 (s) } = \int _{0 } ^\infty \Lambda (x) e^{-sx} 
    \, dx \qquad \text{ for } s\in \bC , \mathrm{Re} \, s > \alpha   \, , 
  \end{equation}
for some  $\alpha <0 $. 
\end{tm}

These equivalences seem to be new. Schoenberg's name is not even
mentioned in ~\cite{Br17-1,Br17-2} on equivalents of the Riemann
hypothesis.

It is interesting  that the characterization of Theorem~\ref{tm:equi2} is
``orthogonal'' to most research on $\zeta $ and to  the  well-known criteria for the Riemann
hypothesis.  Theorem~\ref{tm:equi2} requires  only \emph{the values of
$\zeta $ on the real line} to probe the  secrets of $\zeta $  in the critical strip.
This fact is remarkable, but 
the price to pay  is the added  difficulty  to extract any 
meaningful statements  about  $\xi $ on the critical strip from its
restriction to $\bR $. This seems much harder, if not impossible. 

To work with Theorem~\ref{tm:equi2}, one would need a viable
expression for the Fourier-Laplace transform of $1/\xi $, but there
seems to be none.  The $1$-positivity in \eqref{eq:1} says that
$\Lambda \geq 0$, which is equivalent to the Fourier transform
$\hat{\Lambda } = 1/\xi $ to be positive definite by Bochner's
theorem. Explicitly, we would need to know that,  for all choices of
$c_j \in \bC, \tau _j \in \bR , j=1, \dots , n$, and   all $n\in \bN $, we have $\sum
_{j,k=1} ^n c_j\overline{c_k} \xi (\tfrac{1}{2} + \tau _j -
  \tau _k )\inv  \geq 0$. Not even this property of $1/\xi $ seems to be
known.  It is therefore unlikely that much is  gained by 
Theorems~\ref{tm:equi1} -- \ref{tm:ones}.

By contrast, the Fourier transform of $\Xi (x) $ on
the critical line (!) was
already known to Riemann (see~\cite[2.16.1]{Tit86}) and is the
starting point of a program to prove the Riemann hypothesis that goes
back to Polya~\cite{Pol26}.  After important work of de Bruijn, Hejhal,
and Newman this line of thought has recently culminated in the resolution of the Newman
conjecture by Rodgers and Tao~\cite{RT20}.

\vspace{3mm}

\emph{Some non-trivial \pff s.}
Perhaps Schoenberg  had also the Riemann
hypothesis  in mind, when he investigated  Polya frequency
functions.  
The  examples in~\cite{sch47,sch51} of \tp\ functions smell of  the  zeta function.

(i) The zero set $\{0,-1,-2, \dots \}$
with multiplicity one yields the entire function
\begin{equation}
  \label{eq:n1}
  \Psi (s) = e^{\gamma s} s \prod _{n=1}^\infty (1+\frac{s}{n})\, e^{-s/n}
  \, ,
\end{equation}
where $\gamma $ is the Euler constant. By a classical result $\Psi $ is the reciprocal of the $\Gamma
$-function $\Gamma (s) = \int _0^\infty x^{s-1} e^{-x} \, dx
$. Consequently, the Laplace transform of $\Psi  (s)\inv = \Gamma (s)$
is a \tp\ function. Indeed, using the substitution $x=e^{-t}$ in the
definition of $\Gamma $, one obtains 
\begin{equation}
  \label{eq:n2}
  \Gamma (s) = \int _{-\infty }^\infty e^{-e^{-x}} e^{-sx} \, dx \,
  \qquad \mathrm{Re}\, s >0 \, .
\end{equation}
Theorem~\ref{tm:tp} implies  that
$$
\Lambda (x) = e^{-e^{-x}} 
$$
is \tp . By removing the pole of $\Gamma $ at $0$, we obtain
$$
s\Gamma (s) = \int _{-\infty }^\infty \Lambda'(x)  e^{-sx} \, dx =
\int _{-\infty }^\infty e^{-x} \, e^{-e^{-x}} e^{-sx} \, dx \,  , \qquad
\mathrm{Re} s >-1 \, .
$$
Consequently $\Lambda _1(x) = e^{-x- e^{-x}}$ is a \pff .

(ii) The zero set $\bZ $  with simple
zeros  yields  $\Psi (s) = \frac{\sin \pi s}{\pi}$. By
Theorem~\ref{tm:tp},  $1/\Psi $ is the Laplace transform of a \tp\ function
on a suitable strip of convergence. Schoenberg's calculation yields
the \tp\ function 
$$
\Lambda (x) = \frac{1}{1+e^{-x}} \, .
$$

(iii) Finally the zero set  $\{-n^2: n\in \bN \}$ yields the entire
function 
$$\Psi (s) = s \prod _{n=1}^\infty (1+\frac{s}{n^2}) = -\frac{1}{\pi}
\sqrt{-s} \sin \pi \sqrt{-s} \, . 
$$
The associated \tp\ function is the  Jacobi theta function 
$$
\Lambda (x) =
\begin{cases}
  \sum_{j=-\infty } ^\infty (-1)^j e^{-j^2 x}  & \text{ for } x>0 \\
  0 & \text{ for } x\leq 0 \, .
\end{cases}
$$
 All three functions show up prominently  in the treatment of the
functional equation of the zeta function:  $\Gamma $ is contained  in the
definition of the xi-function, $\sin $ in the formulation of the
functional equation, and a  Jacobi theta function is used in several
proofs of the functional equation  (Riemann's original proof, see~\cite{Tit86}).


\vspace{3mm}

\emph{Intrinsic characterization of \pff s.} The fundamental property
of \pff s is their
smoothing property or \emph{variation diminishing property}. The
relevance of smoothing properties for many applications is outlined in
Schoenberg's survey~\cite{Sch53}. In
this context the variation of a real-valued function on $\bR $
is measured
either by the number of sign changes or by the number of \emph{real}
zeros. Formally, given $f: \bR \to \bR $ let
\begin{equation}
  \label{eq:v1}
  v(f) = \max \# \{n\in \bN : \exists x_j \in \bR,  x_0 < x_1 < \dots < x_n \text{
    with }  f(x_j)f(x_{j+1}) <0\} \, , 
\end{equation}
and let $N(f)$ be the number of \emph{real} zeros of $f$ counted with
multiplicity.

Given a function $\Lambda $, let $T_\Lambda$ be the convolution
operator $T_\Lambda f = f \ast \Lambda $.
Schoenberg's second characterization of \pff s is as
follows~\cite{Sch50}.
\begin{tm} \label{tm:vd}
  Let $\Lambda $ be integrable and continuous. Then  $\Lambda $ is variation diminishing, i.e.,
  $$
  v(T_\Lambda f) \leq v(f) \
  $$
  for all functions that are locally Riemann integrable, \fif\
  either $\Lambda $ or $-\Lambda $  is a
  \pff .
\end{tm}
This characterization is ``intrinsic'' in the sense that it uses only
the properties of the matrices occurring in the definition
\eqref{eq:1} of total positivity.

With a perturbation argument one can replace sign changes with zeros
and obtains the following consequence.
\begin{cor}\label{zdim}
  Let $\Lambda $ be a \pff . Then for every real-valued polynomial $p$
  the convolution $T_\Lambda$ is zero-decreasing, i.e.,
  $$
  N(T_\Lambda p) \leq N(p) \, .
  $$
\end{cor}

\vspace{3mm}

\emph{Intrinsic characterizations of the \lpc .} There are several
characterizations of the \lpc\ that require only their properties as
entire functions. This is part of classical complex analysis and the
results are due to Polya and Schur~\cite{Pol15,PS1914} building on
work of Laguerre, Hadamard, and many others. These results relate the
properties of the zero set to properties of the power series expansion
of an entire function. Before formulating a sequence of equivalences,
we note that every formal power series $F(s) \sim \sum _{j=0}^\infty
a_j s^j$ yields a  differential operator $F(D)p(x) = \sum _{j=0}^\infty
a_j D^jp(x)$ with $D=\frac{d}{dx}$. The differential operator is
well-defined at least  on polynomials, and the mapping $F \mapsto F (D)$ is an
algebra homomorphism and thus  provides  a simple functional calculus.

\begin{tm} \label{tm:lpcc}
Let $\Psi (s) = \sum _{j=0}^\infty \frac{\beta _j}{j!} s^j $ be an
entire function. Then the following are equivalent:

(i) $\Psi $ belongs to the \lpc . 

(ii) $\Psi $ can be approximated uniformly on compact sets by
polynomials with only real zeros. 

(iii) For all $n\in \bN $ the polynomials $p_n (x) = \sum _{j=0}^n \beta _j \binom{n}{j}
x^j$ and $q_n (x) = \sum _{j=0}^n \beta _j \binom{n}{j}
x^{n-j}$  have only real zeros. 

(iv) If $p(x) = \sum _{j=0}^m c_j x^j$ is a polynomial with only real,
\emph{non-positive}  
zeros, then the polynomial $q(x) = \sum \beta _j c_j x^j$ has only
real zeros. 

If, in addition, $\Psi (0)>0$  and $\frac{1}{\Psi (s)} = \sum
_{j=0}^\infty \frac{\gamma _j}{j!} s^j $, then the following property
is equivalent to (i) -- (iv).

(v)  The transform $p \mapsto \frac{1}{\Psi(D)}p$ is zero-decreasing,
i.e., the polynomial $q(x) = \frac{1}{\Psi(D)}p(x ) = \sum
_{j=0}^\infty \frac{\gamma _j}{j!} p^{(j)}(x)$ has at most as many  real zeros
as  $p$ (real-valued):
$$
N\Big(\frac{1}{\Psi (D)}p\Big) \leq N(p) \, .
$$
\end{tm}

Applying condition (iv) to the polynomials $x^{n-1}(1+x)^2$, one
obtains a necessary condition on the Taylor coefficients of a function
in the \lpc , namely the so-called Turan inequalities.

\begin{cor} \label{turan}
  If $\Psi (s) = \sum _{j=0}^\infty \frac{\beta _j}{j!} s^j $ belongs
  to the \lpc , then
  $$
  \beta _n^2 - \beta _{n-1} \beta _{n+1} \geq 0  \qquad \text{ for all }
  n\in \bN
  $$
\end{cor}

Applying condition (v) to polynomials of the form $p(x) =
(\sum_{k=1}^n a_k x^k)^2$ and working out $\Psi   (D)\inv p$, one
obtains the following necessary condition for the \lpc ~\cite[p.~235]{Pol15}.
\begin{cor} \label{hankel}
Assume that  $\Psi $ belongs to the \lpc , $\Psi (0)>0$, $\Psi (s) \neq
e^{as+b}$ and $1/\Psi $ has the Taylor expansion $\frac{1}{\Psi (s)} = \sum
_{j=0}^\infty \frac{\gamma _j}{j!} s^j $. Then for every $n\in \bN $
the $n\times n$ Hankel matrix $\big( \gamma _{j+k} \big) _{j,k =
    0, \dots , n-1} $ is positive definite (and thus invertible). 
\end{cor}
However, the positivity of the Hankel matrices is not sufficient for
$\Psi $ to be in the \lpc , as was proved already by Hamburger~\cite{Ham20}.

Theorem~\ref{tm:lpcc} and its corollaries are all contained in the
seminal papers of Polya and Schur~\cite{Pol15,PS1914} from 1914 and
1915 and have inspired  a century of exciting mathematics. 
Each of the equivalent conditions in 
Theorem~\ref{tm:lpcc} is a point of departure for the study of the
Riemann hypothesis.

 No
list can do justice to all contributions between 1914 and 2020, so let
us mention only a few directions whose origin is in Polya's
work. Further references and  more detailed history can be found in
the cited articles. 

 Condition  (iii) applied to the Riemann function $\Xi $ yields an
important equivalence of the Riemann hypothesis.
The polynomials in condition (iii) are nowadays called Jensen
polynomials.   In modern language
(iii)  says that ``the  Jensen polynomials for the Riemann 
function $\Xi (s)$ must be hyperbolic''. Significant recent progress
on this equivalence  is reported in~\cite{GORZ}. 

The relations between the Jensen polynomials, the multiplier sequences
of condition (iv), and the Turan inequalities  and their
generalizations have been studied in
depth by Craven, Csordas,  and Varga~\cite{Cso15,CC89,CV90} who found many additional
equivalences to the Riemann hypothesis. A particular highlight is
their proof that $\Xi $, or rather the Taylor coefficients of $\Xi
(\sqrt{s})$ satisfy the Turan inequalities~\cite{CNV86}, thereby resolving a 60
year old conjecture going back to --- Polya.  


Finally let us mention that total positivity enters the investigation
of the \lpc\ in yet another way. A entire function belongs to the restricted
\lpc\ defined by \eqref{eq:ne1}, \fif\ the sequence of its Taylor
coefficients $(a_n)$ is a Polya frequency \emph{sequence}~\cite{AS52}. This means
that  the  infinite upper triangular Toeplitz matrix $A$ with entries $A_{jk} =
a_{k-j}$, if $k\geq j$ and $A_{jk} = 0$, if $k<j$ has only positive
minors. This aspect of total positivity  has been used in~\cite{Kat07,Nut13}
for the investigation of the zeta function.

\vspace{3mm}

\emph{From total positivity to the \lpc .}  By
comparing the two intrinsic characterizations in Theorems~\ref{tm:vd}
and \ref{tm:lpcc} one may guess that the respective conditions on zero
diminishing must play the decisive role in the proof of Theorem~\ref{tm:tp}(i).  To give the gist of this
argument, we cannot do better than repeat Schoenberg's beautiful
argument.

First, since $\Lambda $ is assumed to be a \pff , $\Lambda $ must
decay exponentially~\cite[Lemma~2]{sch51}, therefore its moments of
all orders exist. Let
$$
\mu _n  = \int _{\bR } x^n \Lambda (x) \, dx
$$
be the $n$-th moment. By expanding the exponential $e^{-sx}=\sum
_{j=1}^\infty \frac{(-s)^j}{j!} x^j$ we
express the Laplace transform of $\Lambda $ as a power series
\begin{equation}
  \label{eq:nn2}
   \int _{\bR } e^{-sx} \Lambda (x) \, dx = \sum _{j=0}^\infty
  \frac{(-1)^j}{j!} \mu _j \, s^j \, := F(s)  \, .
\end{equation}
Since $\Lambda \not \equiv 0$ and $\Lambda \geq 0$,  we have $F(0) >0$,
and its reciprocal also possesses 
a power series expansion around $0$ with a positive radius of
convergence
$$
\Psi (s) = \frac{1}{F(s)} = \sum _{j=0}^\infty \frac{\beta _j}{j!}\,  s^j
\, .
$$
Next, we consider the convolution of $\Lambda $ with a polynomial $p$ of
degree $N$ and relate it to the moments of $\Lambda $:
\begin{align*}
q(x) &= (\Lambda \ast p)(x) = \int _{\bR } p(x-t) \Lambda (t) \, dt  =\int
 _{\bR } \Big(\sum _{j=0}^N \frac{(-t)^j}{j!} p^{(j)}(x) \Big) \Lambda (t) \, dt \\  
 &= \sum _{j=0 }^n  \frac{(-1)^j}{j!} \mu _j \, p^{(j)}(x) = F(D)p(x)  \, .
\end{align*}
By Corollary~\ref{zdim} the number of real zeros of $q$ (counted with
multiplicity) does not
exceed the number of real zeros of $p$,
\begin{equation}
  \label{eq:nn3}
N(q) = N(F(D) p)  \leq N(p)\, .  
\end{equation}
Using the  functional calculus, we can invert $F(D)$ and  recover $p$ from $q = \Lambda
\ast p$ via
\begin{align*}
  p(x) = \frac{1}{F(D)}q(x) = \Psi (D) q(x) = \sum _{j=0}^\infty
  \frac{\beta _j}{j!}\,  q^{(j)}(x) \, .
\end{align*}
For   the monomial $q(x) = x^n$ we obtain the polynomial
$$
q_n(x) = \Psi (D) x^n = \sum _{j=0}^n \beta _j \binom{n}{j} x^{n-j} \, 
$$
of degree $n$. 
Since $x^n = F(D) q_n$, \eqref{eq:nn3} implies the count of zeros (with
multiplicities) 
$$
 n=  N(x^n) \leq N(q_n) \leq n   \, .
  $$
For every $n$,  $q_n$ therefore has only real zeros. This is precisely
condition (iii) of Theorem~\ref{tm:lpcc}, and we conclude that  $\Psi $ is in the \lpc .

\vspace{3mm}

\emph{Summary.} Schoenberg's characterization of \tp\ functions
implies  a  condition equivalent to the
Riemann hypothesis.  The characterization is
interesting in itself because it involves only the values of the
Riemann zeta function on the real axis. To the best of our knowledge,
the characterization of the \lpc\ by means of \tp\ functions has not
yet been tested on the Riemann zeta function.


\end{document}